\documentclass[number,sort&compress]{elsarticle}

\usepackage{amssymb,amsmath,amsthm,url}
\usepackage[hidelinks]{hyperref}
\usepackage{cleveref}

\crefname{equation}{}{}
\Crefname{equation}{Equation}{Equations}
\crefname{theorem}{Theorem}{Theorems}
\Crefname{theorem}{Theorem}{Theorems}
\crefname{lemma}{Lemma}{Lemmas}
\Crefname{lemma}{Lemma}{Lemmas}
\crefname{proposition}{Proposition}{Propositions}
\Crefname{proposition}{Proposition}{Propositions}
\crefname{corollary}{Corollary}{Corollaries}
\Crefname{corollary}{Corollary}{Corollaries}
\crefname{conjecture}{Conjecture}{Conjectures}
\Crefname{conjecture}{Conjecture}{Conjectures}
\crefname{section}{Section}{Sections}
\Crefname{section}{Section}{Sections}
\crefname{example}{Example}{Examples}
\Crefname{example}{Example}{Examples}

\newcommand{\RR}{\mathbb{R}}

\newcommand{\F}{\mathcal{F}}

\newtheorem{theorem}{Theorem}[section]
\newtheorem{lemma}[theorem]{Lemma}
\newtheorem{proposition}[theorem]{Proposition}
\newtheorem{corollary}[theorem]{Corollary}
\newtheorem{conjecture}[theorem]{Conjecture}
\newdefinition{example}[theorem]{Example}
\newdefinition{problem}[theorem]{Problem}

\numberwithin{equation}{section}

\newcommand{\doi}[1]{\href{http://dx.doi.org/#1}{\texttt{doi:#1}}}

\title{On positive hypergraphs}

\author{Alexander Sidorenko}
\ead{sidorenko.ny@gmail.com}
\address{R\'{e}nyi Institute, Budapest, Hungary}

\date{\today}

\begin{document}

\begin{abstract}
Camarena, Cs\'{o}ka, Hubai, Lippner, and Lov\'{a}sz 
introduced the notion of positive graphs. 
This notion naturally extends to $r$-uniform hypergraphs. 
In the case when $r$ is odd, we prove that 
a hypergraph is positive if and only if its Levi graph is positive.
As an application, we show that 
the $1$-subdivision of $K_{r,r}$ is not a positive graph when $r$ is odd.
\end{abstract}


\maketitle

\section{Introduction}

For a (hyper-)graph $G$, we denote by $V(G)$ and $E(G)$ 
its vertex-set and edge-set, respectively. 
Let ${\rm v}(G):=|V(G)|$ and ${\rm e}(G):=|E(G)|$. 
Let $\F$ denote the class of all bounded measurable functions 
$f:\: [0,1]^2 \to \RR$. 
Let $\F_r^{sym}$ denote the class of all bounded measurable 
symmetric $r$-variate functions 
$f:\: [0,1]^r \to \RR$. 
In particular, $\F_2^{sym}$ is a subclass of $\F$. 

For $f\in\F_2^{sym}$, 
the \emph{homomorphism density} from a graph $G$ to $f$ is defined as
\begin{align}\label{eq:t_G}
     t_G(f) :=  
  \int \prod_{\{u,v\} \in E(G)}
    f(x_u,x_v) \, d\mu^{{\rm v}(G)} ,
\end{align}
where $\mu$ is the Lebesgue measure on $[0,1]$ (see \cite{Lovasz:2010}). 

Let $\pi$ be an automorphism of $G$ such that $\pi = \pi^{-1}$ 
(that is, $\pi$ is an involution). 
Then $\pi$ defines a partition of $V(G)$ into three subsets, 
$V_0,V_+,V_-$, 
where $V_0 = \{v \in V(G): \pi(v)=v\}$ and $\pi(V_+) = V_-$. 
Following the terminology of \cite{Conlon:2017}, 
we call $\pi$ a \emph{stable involution} 
if there are no edges with both ends in $V_0$, 
nor edges with one end in $V_+$ and the other in $V_-$. 

A graph $G$ is called \emph{positive} if 
$t_G(f) \geq 0$ for all $f\in\F_2^{sym}$. 
It is not hard to see that if $G$ has a stable involution, 
then $G$ is positive. 
Camarena, Cs\'{o}ka, Hubai, Lippner, and Lov\'{a}sz \cite{Camarena:2016} 
conjectured that every positive graph must have a stable involution 
(in the language of~\cite{Camarena:2016}, 
such graphs are called ``symmetric''). 
This is known as the \emph{positive graph conjecture}. 

The $1$-subdivision of $K_{r,r}$, where $r$ is even, 
has a stable involution, and hence, is positive.  
When $r \geq 3$ is odd, the $1$-subdivision of $K_{r,r}$ 
does not have a stable involution 
and is one of the simplest graphs 
for which the positive graph conjecture is still open. 

If $G$ is a connected bipartite graph, 
then its vertices are partitioned into two sets, 
and definition~\cref{eq:t_G} naturally extends 
to asymmetric functions $f\in\F$. 
Our first result is 

\begin{theorem}\label{th:1}
If $G$ is a positive bipartite graph, 
then $t_G(f) \geq 0$ for any asymmetric $f\in\F$. 
\end{theorem}

We define positivity for $r$-\emph{graphs} 
(which are $r$-uniform hypergraphs) 
in the same way as it was done for simple graphs. 
Let $H$ be an $r$-graph (possibly, with multiple edges). 
For $f\in\F_r^{sym}$, we define 
\begin{align*}
     t_H(f) :=  
  \int \prod_{\{v_1,v_2,\ldots,v_r\} \in E(H)}
    f(x_{v_1},x_{v_2},\ldots,x_{v_r}) \, d\mu^{{\rm v}(H)} .
\end{align*}
We say that $H$ is \emph{positive} 
if $t_H(f) \geq 0$ for all $f\in\F_r^{sym}$.

The \emph{Levi graph} (or \emph{incidence graph}) of $r$-graph $H$ is 
a bipartite graph $L(H)$ with the bipartition $(V(H),E(H))$, 
where vertices $v \in V(H)$ and $e \in E(H)$ are adjacent 
if and only if $v$ is incident to $e$ in $H$. 
Hence, the adjacency matrix of $L(H)$ is the incidence matrix of $H$. 
We use \cref{th:1} to prove

\begin{theorem}\label{th:2}
Let $r \geq 3$ be odd. 
Then an $r$-graph is positive if and only if 
its Levi graph is positive. 
\end{theorem}

The requirement for $r$ to be odd in \cref{th:2} is essential. 
Indeed, consider an $r$-graph $H$ 
that consists of a single edge. 
Obviously, $H$ is not positive. 
However, $L(H) = K_{r,1}$ is positive when $r$ is even. 

A \emph{homomorphism} from an $r$-graph $H$ 
to an $r$-graph $H'$ 
is a map from $V(H)$ to $V(H')$ which maps edges to edges. 
An \emph{endomorphism} of $H$ is a homomorphism from $H$ to itself. 
An \emph{automorphism} of $H$ is an endomorphism that is a bijection. 
We call an edge $e \in E(H)$ 
\emph{odd} if for any endomorphism $\pi$ of $H$,
$e$ is the image of exactly one of the edges, 
that is, $|\{e' \in E(H): \pi(e')=e\}|=1$. 
In particular, if every pair of vertices is contained 
in at least one of the edges, 
then every endomorphism of $H$ is an automorphism, 
so each edge is odd.

We now give a couple of examples 
to illustrate new opportunities which \cref{th:2} provides. 

\begin{proposition}\label{th:odd}
An $r$-graph that has an odd edge is not positive.
\end{proposition}

Let $n \geq m > k$. The \emph{set-inclusion graph} $I(n,m,k)$ 
is a bipartite graph whose vertices are subsets of $[n]:=\{1,2,\ldots,n\}$
of sizes $m$ and $k$, 
and edges are pairs $X,Y \subseteq [n]$ 
where $|X|=m$, $\,|Y|=k$, $\,Y \subset X$. 
If $m \geq 2k$, and $r = \binom{m}{k}$, then $I(n,m,k)$ 
is the Levi graph of an $r$-graph 
where for every two vertices there is an edge that contains both of them.

\begin{corollary}
If $\binom{m}{k}$ is odd and $n \geq m \geq 2k$, 
then $I(n,m,k)$ is not positive.
\end{corollary}

Let $r$-graph $H$ be $r$-\emph{partite}, that is, 
there exists a homomorphism from $H$ to a single edge. 
Then $H$ has no odd edges. 
Let $H_r$ denote an $r$-graph with $r^2$ vertices 
which are the nodes of an $r \times r$ grid, 
and with $2r$ edges 
which are formed by the horizontal and vertical lines of that grid. 
It is easy to see that $H_r$ is $r$-partite. 
Observe that $L(H_r)$ is the $1$-subdivision of $K_{r,r}\,$. 

\begin{proposition}\label{th:3}
The grid hypergraph $H_r$ is not positive for odd $r \geq 3$.
\end{proposition}

\begin{corollary}
The $1$-subdivision of $K_{r,r}$ is not positive for odd $r \geq 3$.
\end{corollary}

\section{Proofs}

An $r$-dimensional \emph{symmetric tensor} of size $n$ 
is an array of real values $[a_{i_1,i_2,\ldots,i_r}]$ with 
$i_1,i_2,\ldots,i_r \!\in\! \{1,2,\ldots,n\}$ such that 
$a_{i_1,i_2,\ldots,i_r} \!\!=\! 
a_{i_{\sigma(1)},i_{\sigma(2)},\ldots,i_{\sigma(r)}}$ 
for any permutation $\sigma$ of $1,2,\ldots,r$. 
Every symmetric
tensor of size $n$ and dimension $r$ may be uniquely associated 
with a homogeneous polynomial of degree $r$ in $n$ variables. 
It is a classical result known from the XIX century
that there always exists a \emph{symmetric decomposition}
\begin{equation}\label{eq:sym_dec}
  a_{i_1,i_2,\ldots,i_r} = 
  \sum_{j=1}^N \lambda_j \, b_{i_1 j} \, b_{i_2 j} \cdots b_{i_r j}
  \;\;\;\;\;{\rm for}\;\; i_1,i_2,\ldots,i_r \in \{1,2,\ldots,n\},
\end{equation}
where $\lambda_1,\lambda_2,\ldots,\lambda_r$ are reals 
and $[b_{ij}]$ is a real $n \times N$ matrix. 
Lemma~4.2 of \cite{Comon:2008} 
proves a similar statement for complex tensors,
and this proof can be repeated verbatim for real tensors. 
The basic idea is to translate symmetric tensors in polynomials, 
and then the non-existence of a decomposition for a tensor 
translates in the existence of a nonzero polynomial 
that vanishes everywhere on $\RR^n$, giving a contradiction.

We call a symmetric $r$-variate function $h\in\F_r^{sym}$ 
a \emph{step function} of size $n$ 
if there exists a partition $[0,1] = X_1 \cup X_2 \cup \ldots \cup X_n$ 
with $\mu(X_1)=\mu(X_2)=\cdots=\mu(X_n)=1/n$ 
and a symmetric $r$-dimensional tensor $[a_{i_1,i_2,\ldots,i_r}]$ 
of size $n$ such that $h$ is equal to $a_{i_1,i_2,\ldots,i_r}$ 
everywhere on $X_{i_1} \times X_{i_2} \times \cdots \times X_{i_r}$. 

Let $\pi$ be a homomorphism from $r$-graph $H$ to $r$-graph $G$. 
The \emph{homomorphic image} of $H$ (under $\pi$) 
is the subgraph of $G$ whose vertices (edges) 
are images of all vertices (edges) of $H$. 
Let $w: E(G) \to \RR$ define weights of the edges of $G$. 
Then the \emph{weight} of $\pi$ 
is $W(w,\pi) := \prod_{e \in H} w(\pi(e))$. 
If the sum of $W(w,\pi)$ over all homomorphisms $\pi$ from $H$ to $G$ 
is negative, then $H$ is not positive.
Indeed, we may assume that $V(G)=\{1,2,\ldots,N\}$. 
Define a symmetric tensor $[a_{i_1,i_2,\ldots,i_r}]$ with 
$i_1,i_2,\ldots,i_r \in \{1,2,\ldots,N\}$ as follows. 
If $e=\{i_1,i_2,\ldots,i_r\}$ is an edge of $G$, 
we set $a_{i_1,i_2,\ldots,i_r} = w(e)$,
otherwise set $a_{i_1,i_2,\ldots,i_r}=0$. 
Then for a step-function $h\in\F_r^{sym}$ associated with this tensor, 
we get $t_{H_r}(h) < 0$.

\begin{proof}[\bf{Proof of \cref{th:1}}]
If $G_1,G_2,\ldots,G_k$ are connected components of $G$, 
then $t_G(f) = \prod_{i=1}^k t_{G_i}(f)$. 
Thus, it is sufficient to consider the case when $G$ is connected. 
Consider an asymmetric function $f\in\F$. 
Define its ``transpose'' $f^{\mathsf{T}}$ as 
$f^{\mathsf{T}}(x,y) := f(y,x)$. 
Define symmetric function $g\in\F_2^{sym}$ as follows:
\[
  g(x,y) := 
  \begin{cases} 
             0 & \mbox{if } 0 \leq x,y < 1/2; \\
    f(2x,2y-1) & \mbox{if } 0 \leq x < 1/2 \leq y \leq 1; \\
    f(2y,2x-1) & \mbox{if } 0 \leq y < 1/2 \leq x \leq 1; \\
             0 & \mbox{if } 1/2 \leq x,y \leq 1.
  \end{cases}
\]
As $G$ is positive and $g$ is symmetric, $t_G(g) \geq 0$. 
As $G$ is connected, 
\[
  t_G(g) = 2^{-\rm{v}(G)} (t_G(f) + t_G(f^{\mathsf{T}})) .
\]
Hence, $t_G(f) + t_G(f^{\mathsf{T}}) \geq 0$.

Let $h$ be the tensor product of $f$ and $f^{\mathsf{T}}$, 
that is, a function on $[0,1]^4$ defined by 
$h((x,y),(z,w)) := f(x,z) \, f^{\mathsf{T}}(y,w)$. 
Since there is a measure-preserving bijection 
from $[0,1]^2$ onto $[0,1]$, we may regard $h$ as a function on $[0,1]^2$. 
Then $t_G(h) = t_G(f) \cdot t_G(f^{\mathsf{T}})$. 
As $G$ is positive and $h$ is symmetric, 
we get $t_G(f) \cdot t_G(f^{\mathsf{T}}) \geq 0$.
Since both the sum and the product of $t_G(f)$, $t_G(f^{\mathsf{T}})$ 
are nonnegative, $t_G(f)$ itself must be nonnegative.
\end{proof}

\begin{proof}[\bf{Proof of \cref{th:2}}]
Suppose first that $r$-graph $H$ is positive. 
For $f\in\F_2^{sym}$, define $h\in\F_r^{sym}$ by 
\[
  h(x_1,x_2,\ldots,x_r) := \int f(x_1,y) \, f(x_2,y) \cdots f(x_r,y) \,dy.
\]
Then $t_{L(H)}(f) = t_H(h) \geq 0$, 
so $L(H)$ is also positive. 

Now suppose that $H$ is not positive. 
Then there exists a function $g\in\F_r^{sym}$ such that $t_H(g)<0$. 
Since the functional $t_H(\cdot)$ is continuous in the $L_1$ metric, 
there exists an integer $n$ 
and a step function $h\in\F_r^{sym}$ of size $n$ such that $t_H(h)<0$. 
Let $[a_{i_1,i_2,\ldots,i_r}]$ be a symmetric tensor 
associated with the values of $h$. 
As $r$ is odd, we can get rid of coefficients $\lambda_j$ 
in representation \cref{eq:sym_dec} by rescaling 
$c_{ij} = (N\lambda_j)^{1/r} \, b_{ij}\,$: 
\begin{equation*}
  a_{i_1,i_2,\ldots,i_r} = 
  \frac{1}{N}
  \sum_{j=1}^N c_{i_1 j} \, c_{i_2 j} \cdots c_{i_r j}
  \;\;\;\;\;{\rm for}\;\; i_1,i_2,\ldots,i_r \in \{1,2,\ldots,n\}.
\end{equation*}
Consider a partition $[0,1] = Y_1 \cup Y_2 \cup \ldots \cup Y_N$ 
where $\mu(Y_1)=\mu(Y_2)=\cdots=\mu(Y_N)=1/N$ and 
an asymmetric function $f\in\F$ which is equal to $c_{ij}$ 
everywhere on $X_i \times Y_j$ for $i=1,2,\ldots n$, $\;j=1,2,\ldots,N$. 
Then 
\[
h(x_1,x_2,\ldots,x_r) = \int f(x_1,y) \, f(x_2,y) \cdots f(x_r,y) \,dy,
\]
and $t_{L(H)}(f) = t_H(h) < 0$. 
Therefore, by \cref{th:1}, $L(H)$ is not positive. 
\end{proof}

\begin{proof}[\bf{Proof of \cref{th:odd}}]
Let $e$ be an odd edge in $r$-graph $H$. 
Set $w(e)=-1$ and $w(e')=+1$ for edges $e' \neq e$. 
Then for every endomorphism $\pi$ of $H$, we have $W(w,\pi)=-1$. 
Hence, $H$ is not positive. 
\end{proof}

In order to prove \cref{th:3}, we need a few auxiliary results. 
An $r$-graph is called linear if $|e' \cap e''| \leq 1$ 
for any pair of distinct edges $e',e''$. 
An $r$-graph is called a \emph{triangle} 
if it consists of $3$ edges $e',e'',e'''$ such that 
$|e' \cap e''|=|e' \cap e'''|=|e'' \cap e'''|=1$ 
and $e' \cap e'' \cap e''' = \emptyset$. 

\begin{lemma}\label{th:grid}
Let $r \geq 3$. 
Let $\pi$ be a homomorphism from $H_r$ to a linear $r$-graph $G$. 
Then the image of $H_r$ under $\pi$ is either isomorphic to $H_r$, 
or is a single edge, or contains a triangle.
\end{lemma}

\begin{proof}[\bf{Proof}]
We denote the vertices of $H_r$ by $(i,j)$ where $i,j=1,2,\ldots,r$. 
The edges are rows $R_i \!=\! \{(i,1),(i,2),\ldots,(i,r)\}$ 
and columns $C_j \!=\! \{(1,j),(2,j),\ldots,(r,j)\}$. 
If the image of $H_r$ has $r^2$ vertices, 
then it is isomorphic to $H_r$. 
Thus, we may assume that there are two vertices which have the same image. 
Such two vertices can not belong to the same row or column.
Without limiting generality, we may assume that $\pi(1,2)=\pi(2,1)$. 
Notice that the image of $R_1$ and the image of $C_1$ 
have two common vertices $\pi(1,1)$ and $\pi(1,2)=\pi(2,1)$ 
where $\pi(1,1) \neq \pi(1,2)$. 
Since $G$ is linear, then the images of $R_1$ and $C_1$ must coinside. 
Similarly, the images of $R_2$ and $C_2$ must coinside, too. 
We now may renumber rows $R_i$ with $i \geq 3$ 
and columns $C_j$ with $j \geq 3$ in such a way that 
$\pi(1,i) = \pi(i,1)$ for all $i \geq 3$. 
Then the images of $C_i$ and $R_i$ for $i \geq 3$ 
have two common vertices $\pi(i,i)$ and $\pi(1,i)=\pi(i,1)$ 
where $\pi(1,1) \neq \pi(1,i)$. 
By a similar argument, 
we conclude that the images of $R_i$ and $C_i$ coinside. 

We claim that 
if two rows $R_i$ and $R_j$ ($i \neq j$) have the same image, 
then all rows share the same image, 
so the image of $H_r$ is a single edge. 
Indeed, consider a row $R_k$ where $k \neq i,j$. 
As $R_i,R_j,C_i,C_j$ have the same image, 
the image of $R_k$ shares with them two common vertices, 
$\pi(k,i)$ and $\pi(k,j)$. 
As $G$ is linear, the image of $R_k$ must be the same 
as the image of $R_i$ and $R_j$. 

Now assume that all rows have distinct images. 
Since the image of $C_j$ is the same as the image of $R_j$, 
the images of rows $R_i$ and $R_j$ have $\pi(i,j)$ as a common vertex. 
As $G$ is linear, then $|\pi(R_i) \cap \pi(R_j)|=1$ for $i \neq j$. 
Then the images of $R_1,R_2,R_3$ either form a triangle, 
or $|\pi(R_1) \cap \pi(R_2) \cap \pi(R_3)|=1$. 
To finish the proof, we need to exclude the latter case. 
Suppose, $\pi(1,j_1)=\pi(2,j_2)=\pi(3,j_3)$. 
Then by an argument similar to the one we used before, 
the image of $C_{j_1}$ must be the same as the image of $R_2$,
and simultaneously, 
the image of $C_{j_1}$ must be the same as the image of $R_3$.
Hence rows $R_2$ and $R_3$ have the same image, 
which contradicts our assumption.
\end{proof}

\begin{lemma}[{\cite[Section 1.4]{Furedi:2013}}]\label{th:Furedi}
Let $r\geq 3$. For all sufficiently large $n$, 
there exists a linear $r$-graph $G_{r,n}$ 
with $n$ vertices and at least $n^{3/2}$ edges that 
contains neither a triangle nor $H_r$ 
as a subgraph.
\end{lemma}

\begin{proof}[\bf{Proof of \cref{th:3}}]
Let $G_{r,n}$ be the $r$-graph from \cref{th:Furedi}. 
By \cref{th:grid}, every homomorphism from $H_r$ to $G_{r,n}$ 
maps $H_r$ to a single edge of $G_{r,n}$. 
We construct the ``box product'' 
$r$-graph $G_{r,n} \Box G_{r,n}$ 
whose vertex set is $V(G_{r,n}) \times V(G_{r,n})$. 
It has two types of edges. 
The ``horizontal'' ones are 
$\{(x_1,y),(x_2,y),\ldots,(x_r,y)\}$ 
where $\{x_1,x_2,\ldots,x_r\}$ is an edge in $G_{r,n}$. 
The ``vertical'' edges are of form 
$\{(x,y_1),(x,y_2),\ldots,(x,y_r)\}$ 
where $\{y_1,y_2,\ldots,y_r\}$ is an edge in $G_{r,n}$. 
Then ${\rm e}(G_{r,n} \Box G_{r,n}) = 2n {\rm e}(G_{r,n})$. 

We claim that any homomorphism from $H_r$ to $G_{r,n} \Box G_{r,n}$ 
falls to one of two types. 
A homomorphism of the first type maps $H_r$ 
to a single edge of $G_{r,n} \Box G_{r,n}$ 
(horizontal or vertical). 
A homomorphism of the second type maps the vertices of $H_r$ 
to $r^2$ vertices $(x_i,y_j)$ of $G_{r,n} \Box G_{r,n}$ 
($i,j=1,2,\ldots,r$) 
where $\{x_1,x_2,\ldots,x_r\}$ and $\{y_1,y_2,\ldots,y_r\}$ 
are edges in $G_{r,n}$. 
In a homomorphism of the second type, the images of the edges of $H_r$ 
are $r$ horizontal and $r$ vertical edges of $G_{r,n} \Box G_{r,n}$. 
Indeed, let $\pi$ be a homomorphism from $H_r$ to $G_{r,n} \Box G_{r,n}$. 
Suppose, there are two horizontal edges $e',e''$ in $H_r$ such that
$\pi(e')$ is a horizontal edge in $G_{r,n} \Box G_{r,n}$ 
and $\pi(e'')$ is a vertical edge. 
Notice that $|\pi(e') \cap \pi(e'')| \leq 1$, 
and for any $z' \in \pi(e') - \pi(e'')$ and $z'' \in \pi(e'') - \pi(e')$ 
there is no edge in $G_{r,n} \Box G_{r,n}$ 
that would contain both $z'$ and $z''$. 
In this case, $\pi$ would not be able to map vertical edges of $H_r$ 
to the edges of $G_{r,n} \Box G_{r,n}$. 
Thus, we conclude that parallel edges of $H_r$ must be mapped 
to parallel edges of $G_{r,n} \Box G_{r,n}$.
Suppose $\pi$ maps all edges of $H_r$ 
to horizontal edges of $G_{r,n} \Box G_{r,n}$. 
As $H_r$ is connected, each vertex $v$ of $H_r$ must be mapped 
to vertex $(x_v,y_0)$ with fixed $y_0$. 
Then, by the construction of $G_{r,n}$, 
$\pi$ maps all edges of $H_r$ 
to a single horizontal edge of $G_{r,n} \Box G_{r,n}$. 
Similarly, if all edges of $H_r$ are mapped 
to vertical edges of $G_{r,n} \Box G_{r,n}$, 
then the edges of $H_r$ are mapped to a single vertical edge. 
Hence, if $\pi$ does not belong to the first type, 
it must map all horizontal (vertical) edges of $H_r$ to 
horizontal (vertical) edges of $G_{r,n} \Box G_{r,n}$, or vice versa. 
In this case, it is a homomorphism of a second type. 

The number of homomorphisms of the first type is 
$c_r {\rm e}(G_{r,n} \Box G_{r,n}) = 
2c_r n \cdot {\rm e}(G_{r,n})$. 
The number of homomorphisms of the second type is 
$C_r {\rm e}(G_{r,n})^2$. 
Thus, for sufficiently large $n$, 
the number of homomorphisms of the second type is greater. 
We assign weight $+1$ to each horizontal edge of $G_{r,n} \Box G_{r,n}$ 
and $-1$ to each vertical edge. 
Then the weight of a homomorphism of the first type is equal to $+1$. 
As $r$ is odd, 
the weight of a homomorphism of the second type is equal to $-1$. 
Thus, for sufficiently large $n$, 
the sum of the weights of all homomorphisms is negative.
Therefore, $H_r$ is not positive. 
\end{proof}

\section{Concluding remarks}

The notion of stable involution can be generalized to $r$-graphs. 
Let $\pi$ be an involution of $V(H)$. 
Then $\pi$ defines a partition of $V(G)$ into three subsets, 
$V_0,V_+,V_-$, 
where $V_0 = \{v \in V(G): \pi(v)=v\}$ and $\pi(V_+) = V_-$. 
We call $\pi$ a stable involution 
if it preserves the edges,
and there are no edges which are contained in $V_0$ entirely, 
nor edges which intersect both $V_+$ and $V_-$. 
It is not hard to show that if such an involution exists 
then $H$ is positive. 

\begin{conjecture}\label{th:conj}
An $r$-graph is positive if and only if it has a stable involution.
\end{conjecture}

\Cref{th:2} demonstrates that \cref{th:conj} for odd $r$ 
would follow from the positive graph conjecture. 

\medskip

Conlon and Lee \cite{Conlon:2021} proved that every positive graph 
must contain at least one vertex of even degree. 
It implies that a positive $r$-graph with odd $r$ 
also must have a vertex of even degree.

\section*{Acknowledgments}

The author is grateful 
to Joonkyung Lee and Massimiliano Mella for helpful discussions and
to the two anonymous referees for their valuable suggestions.


\begin{thebibliography}{99}

\bibitem{Camarena:2016}
O. A. Camarena, E. Csóka, T. Hubai, G. Lippner, and L. Lov\'{a}sz. 
\newblock Positive graphs. 
\newblock \emph{European J. Combin.}, 52B:290--301, 2016. 
\doi{10.1016/j.ejc.2015.07.007}.

\bibitem{Comon:2008}
P. Comon, G. Golub, L. Lim, and B. Mourrain. 
\newblock Symmetric tensors and symmetric tensor rank. 
\newblock \emph{SIAM J. Matrix Anal. Appl.}, 30:1254--1279, 2008.
\doi{10.1137/060661569}.

\bibitem{Conlon:2017}
D. Conlon and J. Lee. 
\newblock Finite reflection groups and graph norms. 
\newblock \emph{Adv. Math.}, 315:130--165, 2017.
\doi{10.1016/j.aim.2017.05.009}.

\bibitem{Conlon:2021}
D. Conlon and J. Lee. Private communication. 2021.

\bibitem{Furedi:2013}
Z. F\"{u}redi and M. Ruszink\'{o}. 
\newblock Uniform hypergraphs containing no grids. 
\newblock \emph{Adv. Math.}, 240:302--324, 2013.
\doi{10.1016/j.aim.2013.03.009}.

\bibitem{Lovasz:2010}
L. Lov\'{a}sz. 
\newblock \emph{Large networks and graph limits}, 
volume 60 of \emph{Colloquium Publications}. Amer. Math. Soc., 2012. 
\doi{10.1090/coll/060}.

\end{thebibliography}
\end{document}